\documentclass[a4paper,11pt]{amsart} 
\usepackage{amsfonts,amssymb,eucal,amsmath} 

\newcommand{\proofend}{\hspace*{\fill}
$\square$ 
 \normalsize\medskip}

\begin{document}

\title 
{More on counting acyclic digraphs}
\author 
{\bf Valery A. Liskovets}
\address{V.\,A.\,Liskovets:\quad Institute of Mathematics,
National Academy of Sciences, 220072, Minsk, BELARUS}
\email{liskov@im.bas-net.by}

\date{April 14, 2008}  


\begin{abstract}
In this note we derive enumerative formulas for several types of labelled
acyclic directed graphs by slight modifications of the familiar recursive
formula for simple acyclic digraphs. These considerations are motivated by,
and based upon, recent combinatorial results in geometric topology obtained
by S.\,Choi who established exact correspondences between acyclic digraphs and
so-called small covers over hypercubes and related polytopes. In particular,
we show that the number of equivalence classes of small covers over the
cartesian product of $n$ copies of an $r$-simplex is equal to the number
of acyclic $(2^r\!-\!1)$-multidigraphs of order $n$. Asymptotics follows
easily since the main formula is represented by a simple equation in terms
of special generating functions.
\end{abstract} 

\maketitle

\section{Introduction. Simple acyclic digraphs}

A digraph means a finite directed graph without loops and multiple arcs.
Let $a_n$ denote the number of {\sl acyclic} digraphs (that is, digraphs without
oriented cycles) with $n$ labelled vertices, and $a_{n,m}$ be the same for
digraphs with $m$ arcs. Introduce the generating function by $m$
$$
 A_n(x)=\sum\limits_{m=0}^{\binom{n}{2}}a_{n,m}x^m,
$$
where $x$ is a variable. So that, $a_n=A_n(x)\bigl|_{x:=1}=A_n(1)$.
By definition, $a_0=A_0(x)=1.$ The following result obtained
by R.\,Robinson~\cite{Ro73} and R.\,Stanley~\cite{St73}, and rederived
later repeatedly, is well known and has diverse applications:
$$
 a_n=\sum\limits_{t=0}^{n-1}{\binom{n}{t}}(-1)^{n-t-1}2^{t(n-t)}a_t,
\qquad\qquad\quad n\geq 1. \eqno(1)
$$
More generally,
$$
 A_n(x)=\sum\limits_{t=0}^{n-1}{\binom{n}{t}}(-1)^{n-t-1}(1+x)^{t(n-t)}A_t(x),
\qquad n\geq 1. \eqno(2)
$$
These formulas can be proved most easily by the inclusion--exclusion method
(cf.~\cite{Li75}) with respect to the property of a given set of $n-t$
vertices to consist of {\sl sources} (that is, vertices with no incoming
arcs), where $t$ runs over 0,1,\dots,$n-1$ since any nonempty acyclic
digraph possesses at least one source (instead, we could reason dually in terms
of ``sinks'').

Introduce the following special generating functions by $n$ (sometimes
called {\sl graphic} since ${2^{\binom{n}{2}}}$ and ${(1+x)^{\binom{n}{2}}}$
are, respectively, the number and the edge enumerator of labelled graphs
of order $n$):
$$
 A(z)=\sum_{n=0}^\infty\frac{a_nz^n}{n!\,2^{\binom{n}{2}}}\quad
{\rm and}\quad A(x,z)=\sum_{n=0}^\infty\frac{A_n(x)z^n}{n!\,(1+x)^{\binom{n}{2}}}
$$
and
$$
 \Psi(z)=\sum_{n=0}^\infty\frac{(-1)^nz^n}{n!\,2^{\binom{n}{2}}}\quad
{\rm and}\quad \Psi(x,z)=\sum_{n=0}^\infty\frac{(-1)^nz^n}{n!\,(1+x)^{\binom{n}{2}}}.
$$
Then equations (1) and (2) are expressed, respectively, as follows:
$$
 A(z)=1/\Psi(z) \eqno(3)
$$
and
$$
 A(x,z)=1/\Psi(x,z). \eqno(4)
$$
Of course, (3) is a particular case of equation (4) since $A(z)=A(1,z)$
and $\Psi(z)=\Psi(1,z).$

\section{A bicolored variation}

Consider now digraphs with the vertices of two colors, red and blue,
and let $b_n$ denote the number of bicolored acyclic digraphs
with $n$ labelled vertices such that {\sl all red vertices are sources}
(there may be blue sources as well), where the number of vertices
of any color is unspecified. Set $b_0:=1.$

\medskip
{\bf Proposition 1.}
$$
 b_n=\sum\limits_{t=0}^n{\binom{n}{t}}2^{t(n-t)}a_t,\qquad n\ge1.\eqno(5)
$$

Indeed, if $t$ denotes the number of blue vertices then to any
acyclic digraph on the blue vertices one can add $n-t$ red sources
together with arbitrary arcs from red to blue vertices. Relabelling
the vertices arbitrarily and summing over all $t$ we get the rhs
of~(5). And this sum gives rise to $b_n$. \proofend

In view of (3), formula (5) is represented in generating functions as follows:
$$
 B(z)=\Psi(-z)/\Psi(z), \eqno(6)
$$
where
$$
 B(z)=\sum_{n=0}^\infty\frac{b_nz^n}{n!\,2^{\binom{n}{2}}}.
$$

The consideration of $b_n$ is motivated by a curious application to
geometric topology: according to Theorem~3.3 of
S.\,Choi~\cite{Ch08} proved with the help of Burnside's lemma,
the number of $(\mathbb Z_2)^n$-equivariant diffeomorphism classes
of small covers over the $n$-dimensional cube $I^n$ is
$\sum\limits_{t=0}^n{\binom{n}{t}}2^{t(n-t)}a_t\cdot
\prod\limits_{i=0}^{n-1}(2^n-2^i)\bigl/(n!\,2^n)$. See definitions
and details in~\cite{Ch08} and in papers cited therein. Now, by
Proposition~1, this number is equal to
$$
 \frac{b_n}{n!\,2^n}\prod_{i=0}^{n-1}(2^n-2^i).\eqno(7)
$$

{\bf Remark 1.} Some promising quantitative interrelations between diffeomorphism
classes of small covers over hypercubes and {\sl unlabelled} acyclic
digraphs established in~\cite{Ch08} require further investigation.

\section{Acyclic multidigraphs and small covers over powers of a simplex}

Let $k$ be a natural number and set $x:=k$ in formula~(2). Now
$A_n(k)$ can be interpreted as the number acyclic $k$-multidigraphs
with $n$ labelled vertices,
where a multidigraph of multiplicity $k,$ or,
for brevity, a $k$-multidigraph is a multidigraph that has no more than
$k$ parallel arcs from a vertex to another vertex (in particular,
1-multidigraphs are merely simple digraphs; in terms of adjacency
matrices, $k$-multidigraphs are matrices with entries $0,1,\dots,k-1$).
Indeed, formula~(2) with
$x=k$ and with $A_n(k)$ interpreted as stated is provable in the same way.
In particular, ${(1+k)^{t(n-t)}}$ is the number of ways to connect
the selected $n-t$ sources with the remaining $t$ vertices\,\footnote
{Replacing $1+k$ by $1+x+\cdots+x^k$ we would get the corresponding
enumerator by the number of arcs.}.
We will use this simple generalization of acyclic digraphs in order to
derive a closed formula for one more topological enumerator concerning
a uniform generalization of hypercubes.

Let $\Delta^r$ denote an $r$-dimensional simplex and
$\mathbb K_n^r=\Delta^r\times\Delta^r\times\cdots\times\Delta^r
=(\Delta^r)^n$ be the $n$-fold cartesian power of $\Delta^r$;
so that $\mathbb K_n^1$ is a usual $n$-cube. This is a convex polytope
of dimension $rn$. According to Theorem~2.2~\cite{Ch08}, $a_n$
is the number of classes of small covers over the $n$-cube up to
Davis--Januszkiewicz equivalence.

\medskip
{\bf Theorem 1.} {\it The number of classes of small covers over
$\mathbb K_n^r$ up to Davis--Januszkiewicz equivalence $h_n^{(r)}$ is
equal to $A_n(x)\bigl|_{x:=2^r-1}=A_n(2^r-1)$.}

\proof
By Theorem~2.4~\cite{Ch08} applied to this particular uniform case of the product,
$$
 h_n^{(r)}=
\sum\limits_{\Gamma\in\,{\mathfrak A_n}}\prod\limits_{v\in V}(2^r-1)^{d^+_\Gamma(v)},
$$
where $|V|=n,$ $\mathfrak A_n$ is the set of acyclic digraphs on $V$ and
$d^+_\Gamma(v)$ is the outdegree of the vertex $v$ in $\Gamma$. Therefore
$$
 h_n^{(r)}
=\sum\limits_{\Gamma\in\,{\mathfrak A_n}}(2^r-1)^{\sum\limits_{v\in V}{d^+_\Gamma(v)}}
=\sum\limits_{\Gamma\in\,{\mathfrak A_n}}(2^r-1)^{m(\Gamma)}
$$
$$
=\sum\limits_m\sum\limits_{\Gamma\in\,{\mathfrak A_{n,m}}}(2^r-1)^m
=\sum\limits_m a_{n,m}(2^r-1)^m=A_n(2^r-1),
$$
where $m(\Gamma)$ is the number of arcs in $\Gamma$ and $\mathfrak A_{n,m}$
is the set of acyclic digraphs on $V$ with $m$ arcs. \proofend
\medskip

Table~1 contains initial numerical values.
\tabcolsep=0.82ex
\begin{table}[htb!]
\begin{center}
\caption{}
\vspace{-1ex} 
\begin{tabular}{|l|l||r|r|r|r|r|r|}
\hline
$r$&Function    &$n=1$& 2&   3&      4&          5&6\\
\hline
1&$a_n\,\, =\,A_n(1)$&1& 3&  25&    543&      29281&3781503\\
\hline
2&$h_n^{(2)}=A_n(3)$&1& 7& 289&  63487&   69711361&367404658687\\
\hline
3&$h_n^{(3)}=A_n(7)$&1&15&2689&5140479&98267258881&18033699790913535\\
\hline
4&$h_n^{(4)}=A_n(15)$&1&31&23041&{\small365330431}&{\small115851037900801}&{\small705367139018659069951}\\
\hline
 &$b_n$              &2& 8&   74&     1664&          90722&11756288\\
\hline
\end{tabular}
\end{center}
\end{table}

{\bf Remark 2.} In a similar way, one can obtain closed enumerative
formulas for wider classes (including the general class) of cartesian products
of simplices, what requires first to count (say, by the multi-parametric
inclusion--exclusion method) acyclic digraphs with multicolored vertices of
the corresponding specification and with no restrictions on arcs between vertices of the
same or different colors.

\section{Asymptotics}

It is well known and easily provable due to equation~(3) (see,
e.g.,~\cite{Ro73,St73,Li75,BRRW86}) that
$$
 a_n\sim\lambda\,{n!\,2^{\binom{n}{2}}}\omega^{-n}\qquad{\rm as}\quad n\to\infty,
\eqno(8)
$$
where
$$
 \omega=\omega_1=1.4880785455997102947\dots
$$
is the least root of the equation
$$ 
 \Psi(z)=0,
$$
and\, $\lambda=\lambda_1=-1/(\omega\Psi'(\omega))=1.7410611252932298403\dots$
(all calculations were done with the help of Maple). All roots of this equation
are real positive and distinct. Moreover, as pointed out in~\cite{BRRW86},
this is valid in general for the equation
$$
 \Psi(k,z)=0 \eqno(9)
$$
for any fixed $k$ such that $|1+k|>1$. Let $\omega_k$ denote the
least root of~{\rm (9)} and
$$
 \lambda_k:=-{1}\bigl/{\bigl(z\Psi'(k,z)\bigr)}\,\bigl|_{z:=\omega_k}.
$$
Note incidentally that\, $\Psi'(k,z)=-\Psi\bigl(k,\frac{z}{1+k}\bigr).$
Here we need all of this only for $k=2^r-1,$\, $r=1,2,\dots$.
Thus, as an immediate corollary of formulas~(6) and~(4) and
Theorem~1 we obtain (just as in the above-mentioned papers) the
following result.
\medskip

{\bf Proposition 2.} {\it The following asymptotic estimates are valid:
$$
 b_n\sim\Psi(-\omega)\,a_n \qquad {\rm as}\quad n\to\infty,\eqno(10)
$$
where\, $\Psi(-\omega)=3.1135745244678549301\dots$, and for all $r\ge1$,}
$$
 h_n^{(r)}
 =A_n(2^r-1)\sim\lambda_{2^r-1}\frac{n!\,2^{r\binom{n}{2}}}{(\omega_{2^r-1})^n}
  \qquad {\rm as}\quad n\to\infty.\eqno(11)
$$
\proofend

Clearly, $1<\omega_{2^r-1}<1.5$ for all $r$ and $\omega_{2^r-1}$ decreases
to 1 monotonically. Here are numerical values for $r=2,3,4,5$:
$$\begin{array}{ll}
  \omega_3=1.1657706116147275128\dots,\qquad &
 \lambda_3=1.1928652399365987835\dots, \\
  \omega_7=1.0713348333900423361\dots,\qquad &
 \lambda_7=1.0763509327694490247\dots,  \\
  \omega_{15}\!=\!1.0333224614072573348\dots,\qquad &
 \lambda_{15}\!=\!1.0344230890647444796\dots, \\
  \omega_{31}\!=\!1.0161277190328587378\dots,\qquad &
 \lambda_{31}\!=\!1.0163865733813064088\dots .
\end{array}
$$

These constants, as well as $\Psi(-\omega),$ seem to have not
appeared in the literature formerly.


\begin{thebibliography}{9999}
\bibitem
{BRRW86} E.\,A.\,Bender, L.\,B.\,Richmond, R.\,W.\,Robinson and
N.\,C.\,Wormald, The asymptotic number of acyclic digraphs. I.
{\it Combinatorica} {\bf 6}, No.1 (1986), 15--22. MR0856639 (87m:05102),
Zbl\,{\bf 0601}.05025.
\bibitem
{Ch08} Suyoung Choi, The number of small covers over cubes,
Preprint ArXiv:0802.1982v1 [math.GT], 2008.
\bibitem
{Li75} V.\,A.\,Liskovets, On the number of maximal vertices
of a random acyclic digraph, {\it Teor. Veroyatn. Primenen.}
{\bf 20}, No.2 (1975), 412--421 (in Russian; Engl. transl.:
{\it Theory Probab. Appl.} {\bf 20}, No.2 (1975), 401--409).
MR0380925 (52\#1822), Zbl\,{\bf 0362}.60030.
\bibitem
{Ro73} R.\,W.\,Robinson, Counting labeled acyclic digraphs,
{\it New directions in the theory of graphs}, Proc. Third Ann Arbor Conf.
on Graph Theory, 1971 (ed. F.\,Harary), Academic Press, New York (1973),
239--273. MR0363994 (51\#249), Zbl\,{\bf 0259}.05116.
\bibitem
{St73} R.\,P.\,Stanley, Acyclic orientations of graphs, {\it Discrete Math.}
{\bf 5}, No.2 (1973), 171--178. MR0317988 (47\#6537), Zbl\,{\bf 0258}.05113.
\end{thebibliography}
\end{document}